\documentclass[12pt]{article}

\usepackage{natbib}            
\usepackage{graphicx}          
\usepackage{amssymb}                               
\usepackage{amsfonts}

\newcommand{\be}{\begin{equation}}
\newcommand{\ee}{\end{equation}}
\newcommand{\bd}{\begin{displaymath}}
\newcommand{\ed}{\end{displaymath}}
\newcommand{\bea}{\begin{eqnarray}}
\newcommand{\eea}{\end{eqnarray}}

\newcommand{\Hi} {{\cal H}_\infty}

\def\matlab{{\sc \mbox{matlab}}}
\def\hifoo{{\sc \mbox{hifoo}}}
\def\hanso{{\sc \mbox{hanso}}}
\def\mosek{{\sc \mbox{mosek}}}

\begin{document}

\begin{center}
\Large\bf
Multiobjective Robust Control with HIFOO 2.0$^1$
\end{center}

\begin{center}
Suat Gumussoy$^2$, Didier Henrion$^3$, Marc Millstone$^4$, Michael L. Overton$^5$
\end{center}

\footnotetext[1]{
The research of D.\ Henrion was partly supported by project
MSM6840770038 of the Ministry of Education of the
Czech Republic.  The work of M.\ Millstone and M.L.\ Overton was partly
supported by the U.S.\ National Science Foundation under grant DMS-0714321;
views expressed in the paper are those of the authors and not of the NSF.
The work of M.L.\  Overton was also partly funded by Universit\'e Paul Sabatier,
Toulouse, France.
}

\footnotetext[2]{
Katholieke Universiteit Leuven, Department of Computer Science,
Belgium. {\tt suat.gumussoy@cs.kuleuven.be}}

\footnotetext[3]{
LAAS-CNRS, University of Toulouse, France, and
Faculty of Electrical Engineering, Czech Technical University in Prague, Czech
Republic. {\tt henrion@laas.fr}}

\footnotetext[4]{
Courant Institute of Mathematical Sciences,
New York University, USA. {\tt millstone@cims.nyu.edu}}

\footnotetext[5]{
Courant Institute of Mathematical Sciences,
New York University, USA. {\tt overton@cims.nyu.edu}}

\begin{center}
\end{center}

\begin{abstract}
Multiobjective control design is known to be a difficult problem
both in theory and practice.  Our approach is to search for locally
optimal solutions of a nonsmooth optimization problem that is built
to incorporate minimization objectives and constraints for multiple
plants. We report on the success of this approach using our
public-domain \matlab\ toolbox \hifoo\ 2.0, comparing our results with
benchmarks in the literature.
\end{abstract}


\section{Introduction}

Multiobjective control aims at designing a feedback control law
meeting potentially conflicting specifications defined on various
input/output channels.

In the context of linear systems, a standard approach to
multiobjective control is the Lyapunov shaping paradigm proposed in
the mid 1990s by \cite{scherer97}, as an outgrowth of the LMI
(linear matrix inequality) formalism of \cite{boyd94}. Within this
scope, multiobjective controller design boils down to semidefinite
programming (linear programming over the cone of positive
semidefinite matrices) provided all the closed-loop specifications
are certified simultaneously by a unique quadratic Lyapunov
function.
Moreover, the controller is retrieved a posteriori via tedious linear
algebra, and its order is equal to the order of the open-loop plant
plus the order of the weighting functions, which can be quite high
in practice, in contradiction with simplicity of implementability
requirements of embedded control laws.
Another computational approach to multiobjective control exploits
the parametrization of all stabilizing controllers described e.g.
by \cite{sagar85}.
Linear programming can be used in this context to design
controllers, see e.g. \cite{boyd91}, but they are typically of
very high order.

Following a decade of research efforts, these restrictions have been
gradually relaxed (distinct Lyapunov functions for distinct
performance channels, parameter-dependent Lyapunov functions,
decoupling between Lyapunov and controller variables, lower-order
controller design) at the price
of an increased computational burden. The ROMULOC (robust
multiobjective control) toolbox is a recent public-domain \matlab\
implementation of these techniques, see \cite{romuloc06}.

Particular cases of multiobjective robust control problems
include strong stabilization (where a plant must be
stabilized by a controller which is itself stable, see \cite{sagar85})
or simultaneous stabilization (where a single controller
must be found that stabilizes several plants, see \cite{blondel94}).
Most of the algorithms or heuristics available to solve these
problems also typically result in very high order controllers.

In this paper, we introduce the new release 2.0 of our freely
available
package \hifoo, which is aimed at removing the above mentioned limitations
in the context of multiobjective controller design. First, the
controller order is fixed at the outset, allowing for lower-order
controller design. Second, no Lyapunov or lifting variables are
introduced to deal with the conflicting specifications. The
resulting optimization problem is formulated on the controller
coefficients only, resulting in a typically small-dimensional,
nonsmooth, nonconvex optimization problem that does not require the solution of
any large convex subproblems, relieving the computational burden
typical of Lyapunov LMI techniques. Because finding the global
minimum of this optimization problem may be hard, we use an
algorithm that searches only for local minima. While no guarantee
can be made about the behaviour of this algorithm, in practice it
is often possible to determine an acceptable controller quite
efficiently.

See also \cite{anr08} for nonsmooth nonconvex optimization techniques
applied to multiobjective robust control. As far as we know,
no software implementation of these techniques is publicly
available at present.

\section{Problem Formulation}
\label{sec:probform}

The $i$th generalized plant $P^i=(A^i,B^i,C^i,D^i)$ describes the
state-space equations
\begin{eqnarray*}
\label{eq:GeneralizedPlant}
\nonumber \dot{x}^i(t)&=&A^ix^i(t)+B^i_1w^i(t)+B^i_2u^i(t), \\
\nonumber z^i(t)&=&C^i_1x^i(t)+D^i_{11}w^i(t)+D^i_{12}u^i(t), \\
y^i(t)&=&C^i_2x^i(t)+D^i_{21}w^i(t)+D^i_{22}u^i(t), \end{eqnarray*}
where $A^i\in{\mathbb R}^{n^i\times n^i}$,
$D^i_{12}\in{\mathbb R}^{p^i_1\times m_2}$,
$D^i_{21}\in{\mathbb R}^{p_2\times m^i_1}$, with other matrices
having compatible dimensions. The signals $(z^i,w^i,y^i,u^i)$
respectively represent the regulated outputs, the exogenous inputs
(including disturbance and commands), the measured (or sensor)
inputs, and the control inputs.  Let $N$ be the number of plants.

The problem is to choose a \emph{single controller}
\[
K=(A_K,B_K,C_K,D_K)
\]
with state-space equations
\begin{eqnarray*}\label{eq:Controller}
\nonumber \dot{x}_K(t)&=&A_Kx_K(t)+B_Ky(t), \\
u(t)&=&C_Kx_K(t)+D_Ky(t), \end{eqnarray*} where
$A_K\in{\mathbb R}^{n_K\times n_K}$, with $B_K, C_K, D_K$ having
dimensions that are compatible with $A_K$ and the generalized plant
matrices. The controller order $n_K$ is \emph{fixed}, so it can be
specified by the designer. The $\Hi$ \emph{norm} of the $i$th closed loop
system is the norm of the transfer function from input $w^i$ to
output $z^i$; see \cite{zhou96} for details. The \emph{complex
stability radius} is a useful alternative measure when
no $w^i$ and $z^i$ performance channels are specified
and hence the $\Hi$ norm is not defined:
for a stable closed loop system, this is the
largest 2-norm perturbation to the closed loop system
that can be tolerated while guaranteeing that the perturbed system
remains stable.
The \emph{spectral abscissa} of a closed loop system is the largest of
the real parts of its poles (eigenvalues).

Let $\beta_j, j=1,...,N$ each be a real number or $\infty$, and
consider the following optimization problem:
\begin{eqnarray}
        \min_K && \max_{j=1,\ldots,N} \{g_j(K):  \beta_j = \infty\}\label{objfun}\\
        \mathrm{subject~to} &&  g_j(K) \leq \beta_j,  \quad
        j=1,2,\ldots,N,\label{constraints}
\end{eqnarray}
where each $g_j$ is one of the following supported functions
of the closed-loop system for $P^j$,
abbreviated by a single letter as follows:
\begin{itemize}
\item {\tt 'h'}: $\Hi$ norm ($\infty$ if unstable);
\item {\tt 'r'}: reciprocal of complex stability radius ($\infty$ if unstable);
\item {\tt 's'}: spectral abscissa.  
\end{itemize}
The functions $g_j$ are all nonconvex, nonsmooth functions of the
controller matrices. Thus, the optimization problem is potentially
quite difficult.  We focus on two scenarios.

\noindent {\bf Scenario 1}. All $\beta_j=\infty$, and all $g_j$
are {\tt 's'}. Thus the problem is to minimize the maximum of the real
parts of the closed loop poles of all plants.  This approach is
suitable for simultaneous stabilization, as the goal is to move the poles
as far left in the complex plane as possible (in the minmax sense).
If the final objective value is negative, all closed loop plants are stable.

\noindent {\bf Scenario 2}. All $\beta_j=\infty$, and all $g_j$ are
{\tt 'h'}. The problem is to stabilize all plants and minimize the maximum
of the $\Hi$ norms of the closed loop plants.  The advantage over Scenario 1
is that the goal is to not only stabilize the plants, but also optimize their
$\Hi$ performance (in the minmax sense).  The disadvantage is that
evaluating the $\Hi$ norm repeatedly is more time-consuming than computing
the spectral abssissa repeatedly, and I/O performance channels must be
specified.

These two scenarios are the ones for which benchmarking is done in Section
\ref{benchmarks}.  However, the interface to HIFOO is sufficiently flexible
that there are many other alternative ways to call it which may be of
interest to users.  For example:

\noindent {\bf Alternative 3}. All $\beta_j=\infty$, and all $g_j$ are
{\tt 'r'}. The problem is to stabilize all plants and \emph{maximize} the minimum
of the complex stability radii of the closed loop plants (minimize the maximum of
their reciprocals).  The advantage
is that the complex stability radius is a more robust measure of stability
than the spectral abscissa, so this is an appropriate alternative to Scenario 2
when I/O performance channels are not specified.

\noindent {\bf Alternative 4}.  This is a more specific example.
Suppose $\beta_1=\infty$, $\beta_2=100$,
$g_1$ is {\tt 'h'} and $g_2$ is {\tt 'r'}. The problem is to stabilize both
plants and minimize the $\Hi$ norm of the first closed loop plant
subject to the complex stability radius of the second being at least
$0.01$.

If one wishes to impose restrictions, such as stability, on the
\emph{controller} (so-called \emph{strong stabilization}),
it is possible to do so using the multiple-plant model above
by defining a plant $P^i$ so that the closed loop plant
\emph{is equivalent} to the system described by the controller.
However, we provide a more convenient way to specify controller
stability directly, as described below.
Benchmarks assessing the value of a previous
version of \hifoo\ for strong stabilization appear in \cite{cdc08}.

\section{Optimization method}

\hifoo\ 2.0 searches for local minimizers of
(\ref{objfun})-(\ref{constraints}). The algorithm has two phases. In
each phase the main workhorse is the BFGS optimization algorithm,
which is surprisingly effective for nonconvex, nonsmooth
optimization, see \cite{LO-BFGS}. The user can provide an initial guess
for the desired controller (see below); if this is not provided,
\hifoo\ generates randomly generated initial controllers,
and even when an initial guess is provided, \hifoo\ generates some
additional randomly generated initial controllers in case they
provide better results.

The first phase is \emph{stabilization}: BFGS is used to minimize
the maximum of the spectral abscissae of the closed loop plants for
which $g_j$ is either {\tt 'h'} or {\tt 'r'}.  This process terminates as soon
as a controller is found that stabilizes these plants, thus
providing a starting point for which the objective function for the
second phase is finite.

The second phase is \emph{optimization}: BFGS is used to look for a
local minimizer of the following unconstrained problem:
$$
    \min_K ~~ F(K) + \rho \sum_{j=1}^{N} \max(0, g_j(K)-\beta_j),
$$
where $F$ is the objective function defined in (\ref{objfun}) and
$\rho$ is a positive penalty parameter multiplying the sum of the
constraint violations. If BFGS is unable to find a point for which
the constraint violations are zero, the penalty parameter $\rho$ is
increased and the optimization is repeated as needed (unless $F$ is
identically zero, that is all $\beta_j$ are finite).  Although there
are no guarantees, very often this process is quite effective and
reasonably fast. By default, \hifoo\ invokes the
gradient sampling method of \cite{BLO-GradSamp} after BFGS terminates,
but this is generally more time consuming, and can be avoided
as explained in the next section.

\section{User interface}

\hifoo\ is written in \matlab.  The simple call
\begin{center}
\verb@K = hifoo(P, order)@
\end{center}
looks for a controller \verb@K@ solving
Scenario 2 above: stabilize the plants described by \verb@P@ and
minimize the sum of the $\Hi$ norms of the closed loop plants, using
a controller of the specified order.  Here \verb@P@ is a cell array
of plants specified in any of several formats, typically using the
\verb@ss@ class of the \matlab\ Control System Toolbox. A useful
abbreviation is \verb@P{j}='K'@, which specifies that the closed
loop system for the $j$th plant is the controller itself.  If the
\verb@order@ argument is omitted, the default order 0 is used
(static output feedback).

A more general calling sequence is
\begin{verbatim}
 [K, F, viol] = hifoo(P, order, init, fun, ...
                upperbnd, options)
\end{verbatim}
where \verb@P@ and \verb@order@ are as above, \verb@init@ is an
initial guess for the controller (several formats are supported),
\verb@fun@ is a string specifying the characters defining the
supported functions $g_j$ (see Section \ref{sec:probform}) or
a single character if all $g_j$ are the same (default: \verb@'h'@),
\verb@upperbnd@ is an array specifying the upper bounds $\beta_j$
(default: all bounds set to $\infty$),
and \verb@options@ is a structure with various optional fields, some
of which are described below.  The order of the input arguments
does not matter except that \verb@P@ must be first.
The outputs are, in addition to the
controller $K$, the value of the objective function $F(K)$ in
(\ref{objfun}) and a vector of constraint violations
$\max(0,g_j(K)-\beta_j)$. There is a fourth output argument
\verb@loc@ (``local optimality certificate") that can be requested
if gradient sampling is used in addition to BFGS.

Some of the more useful fields in \verb@options@ are:
\begin{itemize}
\item  \verb@options.cpumax@: requests \hifoo\ to quit when the CPU time in seconds
     exceeds this quantity (default: $\infty$)
\item \verb@options.fast@: 1 to use a fast optimization method only (BFGS),
0 to finish optimization with a slower method (gradient
sampling, which may give a better answer) (default: 0, as long as \verb@quadprog@
is in the path; see below)
 \item \verb@options.prtlevel@: one of 0 (no printing), 1 (minimal printing),
2 or 3 (more verbose) (default: 1)
 \item \verb@options.struct@:  specifies \emph{sparsity structure} to be
 imposed on the controller (see documentation for details)
 \item \verb@options.weightNormK@:  weight for adding a penalty on the size of
 the controller to the objective function, specifically
  $(\|A_K\|^2 + \|B_K\|^2 +  \|C_K\|^2 + \|D_K\|^2)^{1/2}$ (default: 0)
 \item \verb@options.augmentHinf@: weight for adding the reciprocal of the complex
   stability radius to the $\Hi$ norm to avoid closed
  loop plants that are only marginally stabilized: applies to all
  plants for which $g_j$ is {\tt 'h'} (default: 0)
 \end{itemize}

\hifoo\ 2.0 uses the following external software:
\begin{itemize}
\item required: \hanso\ 1.0, a hybrid algorithm for non\-smooth optimization,
freely available from the \hifoo\ web page;
\item required: the \matlab\ Control System Toolbox, for
$\Hi$ norm and complex stability radius computation; this also provides
user-friendly system modeling with the {\tt ss} class;
\item optional: \verb@quadprog@ from either the \matlab\ Optimization Toolbox
or \mosek; needed only by the gradient sampling part of the algorithm, which
is not required.
\end{itemize}

\section{Benchmarks}
\label{benchmarks}

We consider a number of simultaneous stabilization problems from the literature;
in each case we wish to stabilize multiple plants with a single controller.
We consider both Scenario 1 (optimizing the spectral abscissa,
by a call such as \verb@K = hifoo(P,'s',n)@, where \verb@n@ is the order
of the controller), and Scenario 2 (optimizing $\Hi$ performance,
by a call such as \verb@K = hifoo(P,'h',n)@, or equivalently \verb@K = hifoo(P,n)@).
Notice that we do not need to explicitly set the upper bounds to $\infty$,
since that is the default value.

The following \matlab\ script illustrates the process, using 
Scenario 1.  There are three first-order plants, given in \cite{jia01}, and
we wish to stabilize them with a first-order controller:
\begin{verbatim}
>> P = {ss(tf([2 -9],[1 -8.8])), ...
        ss(tf([1 2],[1 -6])), ...
        ss(tf([2.5 6],[1 -8]))};
>> K = hifoo(P,'s',1);
...
hifoo: best order 1 controller found
has spectral abscissa -0.284216
...
>> tf(K)
 Transfer function:
s + 1.182
---------
s + 1.595
>> eig(feedback(P{1},-K))
ans =
  -0.2842 + 1.8223i
  -0.2842 - 1.8223i
>> eig(feedback(P{2},-K))
ans =
  1.0e+009 *
   -4.5265
   -0.0000
>> eig(feedback(P{3},-K))
ans =
   -9.7628
   -1.4768
\end{verbatim}
Note that a third-order controller was designed in \cite{jia01}.

The output of \hifoo\ may differ on different runs since the
initialization is done randomly.  However, the output from one run
may be used to initialize a second run on the same problem; the result
cannot be worse.  For example, following the run above with
\begin{verbatim}
>> K = hifoo(P,K,'s',1);
\end{verbatim}
results in
\begin{verbatim}
...
hifoo: best order 1 controller found
has spectral abscissa -0.286524
\end{verbatim}
which is a slight improvement.
Depending on the problem and the initial randomization, several
successive calls to \hifoo\ may be required to obtain a stabilizing
controller.

We have collected various academic and application examples. We now
give benchmark results for Scenarios 1 and~2.

\subsection{Scenario 1: Simultaneous Stabilization}

To evaluate \hifoo, we consider $31$ benchmark problems for
simultaneous stabilization as shown in Table~\ref{table:SimStab}.
There are $11$ problems from applications and $20$
academic test problems. The benchmark problems include 
dynamic, state-feedback and static-output feedback controllers.

\begin{table}[h!]
\begin{center}
\begin{tabular}{@{}l@{}c@{}c@{}c@{}}
  \hline
  Problem Name & $N\times(n^{\max},m_2,p_2)$ & $\:$Known$\:$ & \hifoo\ \\
  \hline
  CRJ$-200$ Aircraft & $6\times(6,1,6)$ & $0$ & $0$ \\
  F$4$E Fighter Aircraft & $4\times(3,1,3)$ & $0$ & $0$ \\
  Gas Turbine Engine & $2\times(10,2,5)$ & $0$ & $0$ \\
  Helicopter Toy & $4\times(3,1,1)$ & $2$ & $1$ \\
  Lane-Keeping of AV & $3\times(2,1,1)$ & $4$ & $1$ \\
  Oblique Wing Aircraft & $64\times(4,1,1)$ & $0$ & $0$ \\
  PFTC & $6\times(4,2,4)$ & $0$ & $0$ \\
  RHM$14$ & $4\times(8,4,4)$ & $26$ & $0$ \\
  Servomotor & $4\times(2,1,1)$ & $1$ & $1$ \\
  Ship-Steering& $2\times(3,1,3)$ & $0$ & $0$ \\
  Stirred-Tank Reactor& $3\times(2,1,1)$ & $0$ & $0$ \\
  \hline
  \hline
  Arehart-Wolovich & $3\times(2,1,1)$ & $1$ & $1$ \\
  Bhattacharyya \textit{et al.}, Ex. 1-1 & $16\times(3,1,1)$ & $0$ & $0$ \\
  Bhattacharyya \textit{et al.}, Ex. 1-2 & $16\times(3,1,1)$ & $1$ & $2$ \\
  Bhattacharyya \textit{et al.}, Ex. 2-1 & $8\times(2,1,1)$ & $0$ & $1$ \\
  Bhattacharyya \textit{et al.}, Ex. 2-2 & $8\times(2,1,1)$ & $1$ & $1$ \\
  Blondel \textit{et al.} & $4\times(1,1,1)$ & $1$ & $0$ \\
  Bredemann, Ex $4.2$ & $3\times(3,1,1)$ & $2$ & $2$ \\
  Bredemann, Ex $5.5$ & $3\times(2,1,1)$ & $3$ & $1$ \\
  Bredemann, Ex $5.6$ & $3\times(2,1,1)$ & $5$ & $0$ \\
  Bredemann, Ex $5.7$ & $3\times(1,1,1)$ & $1$ & $0$ \\
  Cao-Sun & $3\times(2,1,1)$ & $0$ & $0$ \\
  Chen \textit{et al.} & $3\times(2,3,1)$ & $6$ & $0$ \\
  F.-Anaya \textit{et al.}, Ex.$1$ & $80\times(3,1,1)$ & $4$ & $0$ \\
  F.-Anaya \textit{et al.}, Ex.$2$ & $3\times(3,1,1)$ & $6$ & $1$ \\
  F.-Anaya \textit{et al.}, Ex.$3$ & $4\times(4,1,1)$ & $3$ & $1$ \\
  F.-Anaya \textit{et al.}, Ex.$4$ & $5\times(3,1,1)$ & $3$ & $0$ \\
  G\"unde\c{s}-Kabuli & $5\times(10,2,2)$ & $4$ & $1$ \\
  Henrion \textit{et al.}, $1^\textrm{st}$ Ex.  & $3\times(1,1,1)$ & $1$ & $0$ \\
  Henrion \textit{et al.}, $2^\textrm{nd}$ Ex.  & $3\times(1,1,1)$ & $1$ & $0$ \\
  Jia-Ackermann & $3\times(1,1,1)$ & $3$ & $1$ \\
  \hline
  \hline
\end{tabular}
\end{center}
\caption{Benchmarks on simultaneous stabilization: 11 industrial and 20
academic examples. The third and fourth
columns show the lowest stabilizing controller orders in the lit\-era\-ture
and the lowest found by \hifoo.}
\label{table:SimStab}\vspace{-3mm}
\end{table}

Our benchmark results are given in Table~\ref{table:SimStab},
showing the problem name, dimensions (with $n^{\max}=\max\{n^i\}$),
the lowest known controller order from the literature (third column),
and the controller orders obtained by \hifoo\ (final column). 
For each controller order, \hifoo\ was run 10 times;
the best results are reported. Runs for
order $k$ with $k>1$ were initialized with the best controller found
for order $k-1$.  We used \verb@options.fast = 0@ (the default) but all runs
were repeated with \verb@options.fast = 1@, for which the results
were almost the same.

The performance of \hifoo\ is very good considering the large
variety and number of benchmark examples. \hifoo\ successfully
solves the simultaneous stabilization problems with a \emph{low-order}
single controller compared to the existing methods in the
literature. Note that contrary to the other methods in the
literature, \hifoo\ allows the user to set the controller order {\it
a priori}.

In particular, \hifoo\ is very successful in application problems.
\hifoo\ performs better than existing
methods for the application benchmarks Lane-Keeping of Automated
Vehicles and Helicopter Toy. The Rationalized Helicopter Model
(RHM$14$) shows the conservativeness of some of the methods in the
literature, which produce a 26th order stabilizing controller. In
contrast, \hifoo\ stabilizes the same benchmark by a static
controller.

\hifoo\ shows similar performance for academic test problems. For
almost all problems, the results are better than or equivalent to
the results of existing methods. Some of these are Bredemann
Ex.$5.6$, Chen {\it et al.,} F.-Anaya {\it et al.} Ex.$1$ and
Ex.$2$, and G\"unde\c{s}-Kabuli, for all of which the existing methods
stabilize the benchmark problem with a high-order controller and \hifoo\
solves the same problem by a static or first-order controller.

There are two benchmark examples for which \hifoo\ performs slightly
worse than the existing methods, namely Ex.1-1 and 2-1 from
Bhattacharyya et al. which are interval plants. \hifoo\
stabilizes both benchmark problems with a controller whose order is
one more than the order of the controllers in the literature.

\subsection{Scenario 2: Simultaneous $\Hi$ Optimization}

We start with two academic benchmark problems for simultaneous $\Hi$
optimization as shown in Table~\ref{table:SimHinfOpt}. In the first problem, there are
two plants of order two, with optimal full-order $\Hi$ performance (using different
second-order controllers) equal to $1.290$ and $1.245$ respectively. 
\hifoo\ finds a \emph{single} first-order controller with the \emph{same}
$\Hi$ performance (as measured by the maximum of the norms of the two
closed loop plants).  Also, \hifoo\ finds a static controller with performance 1.530.
The previously best known $\Hi$ performance using a single controller
was 1.806, using a fourth-order controller.

The second problem consists of three plants of order two with respective optimal 
$\Hi$ norms $1.290$, $1.245$ and $1.038$ using full-order controllers. 
In the literature a 6th order simultaneously stabilizing controller is known,
 with $\Hi$ performance $1.833$. 
\hifoo\ achieves the optimal $\Hi$ performance with a single first-order controller.

\begin{table}[h!]
\begin{center}
\begin{tabular}{@{}l@{}c@{}c@{}c@{}}
  \hline
  Problem Name & $N\times(n^{\max},m_2,p_2)$ & $\:$Known$\:$ & \hifoo\ \\
  \hline
  Cao-Lam $1$, Ex $2$  & $2\times(2,1,1)$ & $(4,1.806)$ & $(1,1.290)$ \\
  & & & $(0,1.530)$ \\  
  Cao-Lam $2$  & $3\times(2,1,1)$ & $(6,1.833)$ & $(1,1.290)$ \\
  & & & $(0,1.530)$ \\  
  \hline
  \hline
\end{tabular}
\end{center}
\caption{Benchmarks on simultaneous $\Hi$ optimization: two academic examples. 
The third and fourth
columns show the controller order and the minimum value found for the maximum 
of the $\Hi$ norms of the closed loop plants, in the 
lit\-era\-ture and by \hifoo.}
\label{table:SimHinfOpt}\vspace{-3mm}
\end{table}

Table \ref{table:SimHinfOpt2} extends the benchmarks for simultaneous stabilization 
in Table \ref{table:SimStab} to $\Hi$ performance optimization.  
Performance channels are added to each plant in Table~\ref{table:SimStab} using
the \verb@augw@ function in the Control System Toolbox by
\begin{verbatim}
G{k} = augw(P{k},tf(1,[1 1]),[],0.2);
\end{verbatim}
where the weighting functions are $W_1(s)=\frac{1}{s+1}$ and $W_2(s)=0.2$. The dimensions of the plant $P^i$ are 
$A^i\in{\mathbb R}^{n^i+1\times n^i+1}$, $D^i_{12}\in{\mathbb R}^{2p_2\times m_2}$, $D^i_{21}\in{\mathbb R}^{p_2\times m_2}$.

\begin{table}[h!]
\begin{center}
\begin{tabular}{@{}l@{}c@{}c@{}}
  \hline
  \hline
  CRJ$-200$ Aircraft & $(2,2.218)$, $(1,9.333)$, $(0,41.612)$\\  
  F$4$E Fighter Aircraft & $(2,2.993)$, $(1,3.065)$, $(0,6.272)$\\  
  Gas Turbine Engine & $(2,1.000)$, $(1,1.000)$, $(0,1.009)$ \\
  Helicopter Toy & $(3,0.668)$, $(2,0.845)$, $(1,1.079)$ \\
  Lane-Keeping of AV & $(3,0.916)$, $(2,0.947)$, $(1,0.964)$ \\
  Oblique Wing Aircraft & $(2,0.488)$, $(1,0.555)$, $(0,0.959)$ \\
  PFTC & $(2,1.021)$, $(1,1.022)$, $(0,1.027)$ \\
  RHM$14$ & $(2,0.326)$, $(1,0.374)$, $(0,0.867)$ \\
  Servomotor & $(3,0.202)$, $(2,0.202)$, $(1,0.249)$  \\
  Ship-Steering& $(2,1.028)$, $(1,1.033)$, $(0,1.041)$ \\
  Stirred-Tank Reactor& $(2,6.959)$, $(1,6.968)$, $(0,52.921)$  \\
  \hline
  \hline
  Arehart-Wolovich & $(3,2.823)$, $(2,2.823)$, $(1,2.917)$ \\
  Bhattacharyya \textit{et al.}, Ex. 1-1 & $(2,1.000)$, $(1,1.000)$, $(0,1.000)$\\
  Bhattacharyya \textit{et al.}, Ex. 1-2 & $(4,2.236)$, $(3,2.535)$, $(2,16.441)$\\
  Bhattacharyya \textit{et al.}, Ex. 2-1 & $(3,0.202)$, $(2,0.203)$, $(1,0.204)$\\
  Bhattacharyya \textit{et al.}, Ex. 2-2 & $(3,3.725)$, $(2,3.739)$, $(1,3.745)$\\
  Blondel \textit{et al.} & $(2,3.053)$, $(1,3.053)$, $(0,3.053)$ \\
  Bredemann, Ex $4.2$ & $(4,1.623)$, $(3,1.736)$, $(2,1.772)$  \\
  Bredemann, Ex $5.5$ & $(3,1.515)$, $(2,1.522)$, $(1,1.682)$  \\
  Bredemann, Ex $5.6$ & $(2,31.613)$, $(1,31.613)$, $(0,31.613)$ \\
  Bredemann, Ex $5.7$ & $(2,1.073)$, $(1,1.173)$, $(0,7.517)$ \\
  Cao-Sun & $(2,0.201)$, $(1,0.202)$, $(0,0.202)$ \\
  Chen \textit{et al.} &  $(2,0.595)$, $(1,0.875)$, $(0,1.000)$\\
  F.-Anaya \textit{et al.}, Ex.$1$ & $(2,0.200)$, $(1,0.200)$, $(0,0.200)$\\
  F.-Anaya \textit{et al.}, Ex.$2$ & $(3,0.256)$, $(2,0.256)$, $(1,0.256)$\\
  F.-Anaya \textit{et al.}, Ex.$3$ & $(3,1.000)$, $(2,1.000)$, $(1,1.000)$\\
  F.-Anaya \textit{et al.}, Ex.$4$ & $(2,0.201)$, $(1,0.203)$, $(0,0.203)$\\
  G\"unde\c{s}-Kabuli & $(3,1.000)$, $(2,1.000)$, $(1,1.000)$\\
  Henrion \textit{et al.}, $1^\textrm{st}$ Ex. & $(2,0.7155)$, $(1,0.759)$, $(0,1.059)$\\
  Henrion \textit{et al.}, $2^\textrm{nd}$ Ex. & $(2,3.044)$, $(1,3.044)$, $(0,3.455)$\\
  Jia-Ackermann & $(3,2.556)$, $(2,3.080)$, $(1,13.594)$\\
  \hline
  \hline
\end{tabular}
\end{center}
\caption{Benchmarks on simultaneous $\Hi$ optimization: 11 industrial and 20
academic examples. The second
column shows the lowest three stabilizing controller orders found by \hifoo\ 
and the corresponding max\-imum of the $\Hi$ norms of the \mbox{closed loop plants}.}
\label{table:SimHinfOpt2}\vspace{-3mm}
\end{table}

These results clearly demonstrate that \hifoo\ is very effective
over various types of benchmark examples including industrial
application and academic test problems. 

\hifoo\ 2.0 is available under the GNU Public License at
\begin{verbatim}
www.cs.nyu.edu/overton/software/hifoo
\end{verbatim}
Further information on references for the benchmark problems and
other methods will be provided in a technical report that will
be made available at this website.

In conclusion, \hifoo, by
allowing the designer to specify the controller order, is a useful
tool for simultaneous stabilization and simultaneous $\Hi$ optimization.

\end{document}